\documentstyle[amssymb,amsfonts]{amsart}

\def\normo#1{\left\|#1\right\|}

\def\brk#1{\left(#1\right)}

\def\rev#1{\frac{1}{#1}}
\def\half#1{\frac{#1}{2}}
\def\norm#1{\|#1\|}
\def\jb#1{\langle#1\rangle}
\def\wt#1{\widetilde{#1}}
\def\wh#1{\widehat{#1}}

\numberwithin{equation}{section}

\newcommand{\N}{{\mathbb N}}
\newcommand{\T}{{\mathbb T}}

\newcommand{\R}{{\mathbb R}}
\newcommand{\C}{{\mathbb C}}
\newcommand{\Z}{{\mathbb Z}}
\newcommand{\ft}{{\mathcal{F}}}

\newcommand{\les}{{\lesssim}}
\newcommand{\ges}{{\gtrsim}}
\newcommand{\ra}{{\rightarrow}}

\newcommand{\Sch}{{\mathcal{S}}}
\newcommand{\supp}{{\mbox{supp}}}

\theoremstyle{plain}
  \newtheorem{theorem}[subsection]{Theorem}
  \newtheorem{proposition}[subsection]{Proposition}
  \newtheorem{lemma}[subsection]{Lemma}

\theoremstyle{remark}
  \newtheorem{remark}[subsection]{Remark}

\theoremstyle{definition}

\newenvironment{proof}{\noindent {\bf Proof.} }{\endprf\par}
\def \endprf{\hfill  {\vrule height6pt width6pt depth0pt}\medskip}

\begin{document}
\title[Korteweg-de Vries equation]{Global Well-posedness of Korteweg-de Vries equation in $H^{-3/4}(\R)$}

\author{Zihua Guo}
\address{LMAM, School of Mathematical Sciences, Peking University, Beijing
100871, China}

\email{zihuaguo@@math.pku.edu.cn}

\begin{abstract}
We prove that the Korteweg-de Vries initial-value problem is
globally well-posed in $H^{-3/4}(\R)$ and the modified Korteweg-de
Vries initial-value problem is globally well-posed in $H^{1/4}(\R)$.
The new ingredient is that we use directly the contraction principle
to prove local well-posedness for KdV equation in $H^{-3/4}$ by
constructing some special resolution spaces in order to avoid some
'logarithmic divergence' from the high-high interactions. Our local
solution has almost the same properties as those for $H^s\ (s>-3/4)$
solution which enable us to apply the I-method to extend it to a
global solution.
\end{abstract}



\maketitle

\section{Introduction}

This paper is mainly concerned with the well-known open problem for
the Korteweg-de Vries equation: global well-posedness holds in
$H^{-3/4}$ (cf. \cite{I-method}, \cite{CCT})? The Korteweg-de Vries
(KdV) equation
\begin{eqnarray}\label{eq:kdv}
\left\{
\begin{array}{l}
u_t+u_{xxx}-3(u^2)_x=0,\quad u(x,t):\R\times \R\rightarrow \R,\\
u(x,0)=\phi(x)\in H^s(\R),
\end{array}
\right.
\end{eqnarray}
has attracted extensive attentions, since it was first derived by D.
J. Korteweg and G. de Vries \cite{KdV} as a model for unidirectional
propagation of nonlinear dispersive long waves. A large amount of
works have been devoted to the Cauchy problem \eqref{eq:kdv} and
meanwhile many useful tools and methods were discovered and
developed. We only mention here the most recent results concerned
with the well-posedness. A first result by using contraction
principle was due to Kenig, Ponce and Vega \cite{KPVCPAM93} who
obtained local well-posedness in $H^s$ for $s>3/4$. Bourgain
\cite{Bour} extended this result to global well-posedness in $L^2$
by developing the $X^{s,b}$ space. Then by developing the bilinear
estimates in $X^{s,b}$ space Kenig, Ponce and Vega \cite{KPVJAMS96}
were able to prove local well-posedness in $H^{s}$ for $s>-3/4$ and
Colliander, Keel, Staffilani, Takaoka and Tao \cite{I-method}
extended it to a global result where $I-method$ was introduced. It
is worth noting that $H^{-3/4}$ is the sharp regularity for a strong
well-posedness. Christ, Colliander, and Tao \cite{CCT} proved that
the KdV solution map of \eqref{eq:kdv} fails to be uniformly
continuous in $H^s$ for $s<-3/4$ which was first proved by Kenig,
Ponce and Vega \cite{KPVDMJ01} for the complex-valued problem.

On the other hand, analogous analysis above were also devoted to the
modified KdV (mKdV) equation
\begin{eqnarray}\label{eq:mkdv}
u_t+u_{xxx}\pm 2(u^3)_x=0,\quad u(x,0)=\phi(x).
\end{eqnarray}
It is well-known that under the Miura transform
$v=M(u)=\partial_xu+u^2$ for $'-'$ ($\partial_xu+iu^2$ for '+')
where $u$ is a solution to \eqref{eq:mkdv}, then $v$ satisfies the
KdV equation \eqref{eq:kdv} (with some different coefficient in the
nonlinearity). Thus we see from the Miura transform that mKdV
requires 1-order higher regularity than the KdV equation. But mKdV
has better result at the endpoint $s=1/4$, since it has better
high-high interactions. It is known that the mKdV equation
\eqref{eq:mkdv} is locally well-posed in $H^s$ for $s\geq 1/4$
\cite{KPVJAMS96} and globally well-posed for $s>1/4$
\cite{I-method}. The $H^{1/4}$ well-posedness for the defocusing
mKdV equation combined with the Miura transform established a local
well-posedness result for KdV in $H^{-3/4}$ \cite{CCT}. Global
well-posedness of KdV at $H^{-3/4}$ and for mKdV in $H^{1/4}$ remain
open problems.

In this paper we intend to study the global well-posedness of KdV at
$s=-3/4$ and of mKdV at $s=1/4$. Generally, there are two approaches
to this problem. The main reason that the $H^{-3/4}$ local solution
of KdV in \cite{CCT} can not be extended to a global solution by
using the I-method is that the local solution doesn't have as many
properties as the $H^s$ local solution for $s>-3/4$. These
properties are crucial in I-method to estimate the multi-linear
forms and hence control the increase of the modified energy.
However, on the other hand, the $H^{1/4}$ local solution of mKdV was
derived by direct contraction principle and thus has many
properties. Therefore it is expected that one may follow the ideas
of I-method to directly study mKdV to get global well-posedness in
$H^{1/4}$ and then prove global well-posedness for KdV. One need to
estimate a four or perhaps higher linear form.  The other approach
is to prove a strong local well-posedness for the KdV at $s=-3/4$.
This is possible because the solution map is locally uniformly
continuous (actually analytic) in $H^{-3/4}$ although the uniform
continuity fails in $H^s$ if $s<-3/4$. We will adapt the second
approach. We construct a strong local solution by contraction
principle. Now we state our main results:

\begin{theorem}\label{thmlwp} Assume $\phi \in H^{-3/4}$. Then

(a) Existence. There exist $T=T(\norm{\phi}_{H^{-3/4}})>0$ and a
solution $u$ to the Cauchy problem \eqref{eq:kdv} satisfying
\[u\in \bar{F}^s(T) \subset C([-T,T]:H^{-3/4}).\]

(b) Uniqueness. The solution mapping $S_T:\phi\rightarrow u$ is the
unique extension of the classical solution $H^\infty\rightarrow
C([-T,T]:H^\infty)$.

(c) Lipschitz continuity. For any $R>0$, the mapping
$\phi\rightarrow u$ is Lipschitz continuous from $\{u_0\in
H^{-3/4}:\norm{u_0}_{H^{-3/4}}<R\}$ to $C([-T,T]:H^{-3/4})$.

(d) Persistence of Regularity. If in addition $\phi \in H^s$ for
some $s>-3/4$, then the solution $u \in H^s$.
\end{theorem}

We remark that Theorem \ref{thmlwp} also holds for the
complex-valued KdV equation. From Theorem \ref{thmlwp} (b) and (d)
we get that our local solution coincide with the one in \cite{CCT},
and we also prove it belongs to a strong class $\bar{F}^s(T)$. With
this we are able to use I-method to extend it to a global solution.

\begin{theorem}\label{thmgwp}
The KdV equation \eqref{eq:kdv} is globally well-posed in $H^{-3/4}$
and the mKdV  equation \eqref{eq:mkdv} is globally well-posed in
$H^{1/4}$.
\end{theorem}

By 'globally well-posed' for KdV in Theorem \ref{thmgwp} we mean $T$
can be any large real number in Theorem \ref{thmlwp} and similarly
for mKdV. We will prove Theorem \ref{thmgwp} by using the I-method
and a variant version of Theorem \ref{thmlwp}. Now we sketch our
ideas in proving Theorem \ref{thmlwp}. Our starting point is the
fundamental bilinear estimates in $X^{s,b}$ space (see next section
for the definition of $X^{s,b}$ ):
\begin{eqnarray}\label{eq:bilinearXsb}
\norm{\partial_x(uv)}_{X^{s,b-1}}\leq
C\norm{u}_{X^{s,b}}\norm{v}_{X^{s,b}}.
\end{eqnarray}
The bilinear estimates \eqref{eq:bilinearXsb} play a crucial role in
\cite{KPVJAMS96} to apply a fixed point argument. It was proved in
\cite{KPVJAMS96} that \eqref{eq:bilinearXsb} hold for some $b>1/2$
if $s>-3/4$ and fail for any $b$ if $s<-3/4$. For $s=-3/4$, the
bilinear estimates \eqref{eq:bilinearXsb} also fail for any $b$
which was due to Nakanishi, Takaoka and Tsutsumi \cite{NTT}. In view
of this, we study instead the bilinear estimates in $l^1$-type
$X^{s,b}$ space $F^s$ inspired by our work \cite{GW}. We find that
the bilinear estimates in $F^{-3/4}$ almost hold except some
potential logarithmic divergences from $high\times high \rightarrow
low $ interactions. Fortunately, we find that we are still able to
control the $high\times low$ interactions if assuming some weaker
structure on the low frequency. Using the weaker structure we are
then able to control the $high\times high$ interactions.

Our methods can also be adapted to the other similar problems where
some logarithmic divergences appear in the high-high interactions.
One interesting problem is the global well-posedness for the
KdV-burger equations in $H^{-1}$:
\begin{eqnarray}\label{eq:kdv-b}
u_t+u_{xxx}-\epsilon u_{xx}-3(u^2)_x=0,\quad u(x,0)=\phi(x).
\end{eqnarray}
The equation \eqref{eq:kdv-b} is showed in \cite{MRIMRN} to be
globally well-posed in $H^s$ for $s>-1$ and $C^2$ ill-posed for
$s<-1$, but $H^{-1}$ well-posedness remains a open problem. Some new
ideas should be developed to this problem. One may also follow the
methods here and the ideas in \cite{GW} to prove the inviscid limit
in $C([-T,T]:H^{-3/4})$ as $\epsilon$ tends to zero. We do not
pursue this in this paper.

\begin{remark}
After this paper was published, the author learned that Nobu
Kishimoto \cite{Kishimoto} also obtained similar results with
different resolution spaces by using some ideas in \cite{BT}. The
author would like to thank him for valuable discussion.
\end{remark}

The rest of the paper is organized as following. In Section 2 we
present some notations and Banach function spaces. We present some
dyadic bilinear estimates in Section 3. The proof of Theorem
\ref{thmlwp} and \ref{thmgwp} is given in Section 4.

\section{Notation and Definitions} \label{notation}
For $x, y\in \R$, $x\sim y$ means that there exist $C_1, C_2 > 0$
such that $C_1|x|\leq |y| \leq C_2|x|$. For $f\in \Sch'$ we denote
by $\widehat{f}$ or $\ft (f)$ the Fourier transform of $f$ for both
spatial and time variables,
\begin{eqnarray*}
\widehat{f}(\xi, \tau)=\int_{\R^2}e^{-ix \xi}e^{-it \tau}f(x,t)dxdt.
\end{eqnarray*}
We denote  by $\ft_x$ the Fourier transform on spatial variable and
if there is no confusion, we still write $\ft=\ft_x$. Let
$\mathbb{Z}$ and $\mathbb{N}$ be the sets of integers and natural
numbers, respectively. $\Z_+=\N \cup \{0\}$. For $k\in \Z_+$ let
\[{I}_k=\{\xi: |\xi|\in [2^{k-1}, 2^{k+1}]\}, \ k\geq 1; \quad I_0=\{\xi: |\xi|\leq 2\}.\]
Let $\eta_0: \R\rightarrow [0, 1]$ denote an even smooth function
supported in $[-8/5, 8/5]$ and equal to $1$ in $[-5/4, 5/4]$. We
define $\psi(t)=\eta_0(t)$. For $k\in \Z$ let
$\eta_k(\xi)=\eta_0(\xi/2^k)-\eta_0(\xi/2^{k-1})$ if $k\geq 1$ and
$\eta_k(\xi)\equiv 0$ if $k\leq -1$. For $k\in \Z$ let
$\chi_k(\xi)=\eta_0(\xi/2^k)-\eta_0(\xi/2^{k-1})$. Roughly speaking,
$\{\chi_k\}_{k\in \mathbb{Z}}$ is the homogeneous decomposition
function sequence and $\{\eta_k\}_{k\in \mathbb{Z}_+}$ is the
non-homogeneous decomposition function sequence to the frequency
space. For $k\in \Z$ let $P_k$ denote the operator on $L^2(\R)$
defined by
\[
\widehat{P_ku}(\xi)=\eta_k(\xi)\widehat{u}(\xi).
\]
By a slight abuse of notation we also define the operator $P_k$ on
$L^2(\R\times \R)$ by the formula $\ft(P_ku)(\xi,
\tau)=\eta_k(\xi)\ft (u)(\xi, \tau)$. For $l\in \Z$ let
\[
P_{\leq l}=\sum_{k\leq l}P_k, \quad P_{\geq l}=\sum_{k\geq l}P_k.
\]
Thus we see that $P_{\leq 0}=P_0$.

For $u_0\in \Sch'(\R)$, we denote by $W(t)u_0=e^{-t\partial_x^3}u_0$
the free solution of linear Airy equation which is defined as
\[
\ft_x(W(t)\phi)(\xi)=\exp[i\xi^3t]\widehat{\phi}(\xi), \  \forall \
t\in \R.
\]
 We define the Lebesgue spaces $L_{t\in I}^qL_x^p$
and $L_x^pL_{t\in I}^q$ by the norms
\begin{equation}
\norm{f}_{L_{t\in I}^qL_x^p}=\normo{\norm{f}_{L_x^p}}_{L_t^q(I)},
\quad \norm{f}_{L_x^pL_{t\in
I}^q}=\normo{\norm{f}_{L_t^q(I)}}_{L_x^p}.
\end{equation}
If $I=\R$ we simply write $L_{t}^qL_x^p$ and $L_x^pL_{t}^q$. We will
make use of the $X^{s,b}$ norm associated to the KdV equation
\eqref{eq:kdv} which is given by
\begin{eqnarray*}
\norm{u}_{X^{s,b}}=\norm{\jb{\tau-\xi^3}^b\jb{\xi}^s\widehat{u}(\xi,\tau)}_{L^2(\R^2)},
\end{eqnarray*}
where $\jb{\cdot}=(1+|\cdot|^2)^{1/2}$. The spaces $X^{s,b}$ turn
out to be very useful in the study of low-regularity theory for the
dispersive equations. These spaces were first used to systematically
study nonlinear dispersive wave problems by Bourgain [5] and
developed by Kenig, Ponce and Vega \cite{KPVJAMS96} and Tao
\cite{Taokz}. Klainerman and Machedon \cite{KlMa} used similar ideas
in their study of the nonlinear wave equation.

In applications we usually apply $X^{s,b}$ space for $b$ is very
close to $1/2$. In the case $b=1/2$ one has a good substitute-$l^1$
type $X^{s,b}$ space. For $k\in \Z_+$ we define the dyadic
$X^{s,b}$-type normed spaces $X_k=X_k(\R^2)$,
\begin{eqnarray}
X_k=\left\{f\in L^2(\R^2):
\begin{array}{l}
f(\xi,\tau) \mbox{ is supported in } I_k\times\R \mbox{ and }\\
\norm{f}_{X_k}=\sum_{j=0}^\infty
2^{j/2}\norm{\eta_j(\tau-\xi^3)\cdot f}_{L^2}.
\end{array}
\right\}
\end{eqnarray}
Then we define the $l^1$-analogue of $X^{s,b}$ space ${F}^s$ by
\begin{eqnarray}
\norm{u}_{{F}^s}^2=\sum_{k \geq
0}2^{2sk}\norm{\eta_k(\xi)\ft(u)}_{X_k}^2.
\end{eqnarray}
Structures of this kind of spaces were introduced, for instance, in
\cite{Tata}, \cite{IKT} and \cite{In-Ke} for the BO equation. The
space $F^s$ is better than $X^{s,1/2}$ in many situations for
several reasons. $F^s$ can be embedded into $C(\R;H^s)$ and into the
Strichartz-type space, say $L^p_tL^q_x$ as $X^{s,1/2+}$ (see Lemma
\ref{lemmaextension} below). On the other hand, it has the same
scaling in time as $X^{s,1/2}$, which was recently exploited by us
\cite{GW} in the inviscid limit problem for KdV-burgers equation.
This is similar to the difference between $\dot{B}_{2,1}^{n/2}$ and
$\dot{H}^{n/2}$. Tao \cite{Taoscatteringgkdv} used a homogenous
version to prove scattering for generalized KdV equation for small
critical data. From the definition of $X_k$, we see that for any
$l\in \Z_+$ and $f_k\in X_k$ (see also \cite{IKT}),
\begin{eqnarray}
\sum_{j=0}^\infty 2^{j/2}\normo{\eta_j(\tau-\xi^3)\int
|f_k(\xi,\tau')|2^{-l}(1+2^{-l}|\tau-\tau'|)^{-4}d\tau'}_{L^2}\les
\norm{f_k}_{X_k}.
\end{eqnarray}
Hence for any $l\in \Z_+$, $t_0\in \R$, $f_k\in X_k$, and $\gamma
\in \Sch(\R)$, then
\begin{equation}\label{eq:propertyXk2}
\norm{\ft[\gamma(2^l(t-t_0))\cdot \ft^{-1}f_k]}_{X_k}\les
\norm{f_k}_{X_k}.
\end{equation}
In order to avoid some logarithmic divergence, we need to use a
weaker norm for the low frequency
\begin{eqnarray*}
\norm{u}_{\bar{X}_0}=\norm{u}_{L_x^2L_t^\infty}.
\end{eqnarray*}
It is easy to see from Proposition \ref{propXkembedding} that
\begin{eqnarray}\label{eq:xbar0x0}
\norm{\eta_0(t)P_{\leq 0}u}_{\bar{X}_0}\les \norm{P_{\leq
0}u}_{{X}_0}.
\end{eqnarray}
On the other hand, for any $1\leq q\leq \infty$ and $2\leq r\leq
\infty$ we have
\begin{eqnarray}\label{eq:Xbar0Lqr}
\norm{P_{\leq0}u}_{L_{|t|\leq T}^qL_x^r\cap L_x^rL_{|t|\leq
T}^q}\les_T \norm{P_{\leq0}u}_{L_x^2L_{|t|\leq T}^\infty}.
\end{eqnarray}
For $-3/4\leq s\leq 0$, we define the our resolution spaces
\begin{eqnarray*}
\bar{F}^s=\{u\in \Sch'(\R^2):\norm{u}_{\bar{F}^s}^2=\sum_{k \geq
1}2^{2sk}\norm{\eta_k(\xi)\ft(u)}_{X_k}^2+\norm{P_{\leq
0}(u)}_{\bar{X}_0}^2<\infty\}.
\end{eqnarray*}
For $T\geq 0$, we define the time-localized spaces $\bar{F}^{s}(T)$:
\begin{eqnarray}
\norm{u}_{\bar{F}^{s}(T)}=\inf_{w\in \bar{F}^{s}}\{\norm{P_{\leq
0}u}_{L_x^2L_{|t|\leq T}^\infty}+\norm{P_{\geq 1}w}_{\bar{F}^{s}}, \
w(t)=u(t) \mbox{ on } [-T, T]\}.
\end{eqnarray}

Let $a_1, a_2, a_3\in \R$. It will be convenient to define the
quantities $a_{max}\geq a_{med}\geq a_{min}$ to be the maximum,
median, and minimum of $a_1,a_2,a_3$ respectively. Usually we use
$k_1,k_2,k_3$ and $j_1,j_2,j_3$ to denote integers, $N_i=2^{k_i}$
and $L_i=2^{j_i}$ for $i=1,2,3$ to denote dyadic numbers.

\section{Dyadic Bilinear Estimates}

In this section we prove some dyadic bilinear estimates which are
crucial for applying contraction principle in the next section. We
will need the estimates for the free solution to the KdV equation.
We recall in the following lemma the Strichartz estimates, maximal
function estimates and smoothing effect estimate for the free KdV
solution and refer the readers to \cite{GPW,KPVIUMJ91,KPVJAMS91} for
their proof.
\begin{lemma}[Estimates for free KdV solutions]\label{lemmakdvfree}
Let $I\subset \R$ be a interval with $|I|\les 1$ and $k\in \Z_+$.
Then for all $\phi \in \Sch(\R)$ we have
\begin{eqnarray}
&&\norm{W(t)\phi}_{L_t^qL_x^r}\les \norm{\phi}_{L^2},\\
&&\norm{W(t)P_k(\phi)}_{L_x^2L_{t\in I}^\infty}\les 2^{3k/4}
\norm{\phi}_{L^2},\\
&&\norm{W(t)\phi}_{L_x^4L_{t}^\infty}\les
\norm{\phi}_{\dot{H}^{1/4}},\\
&&\norm{W(t)\phi}_{L_x^\infty L_{t}^2}\les
\norm{\phi}_{\dot{H}^{-1}},
\end{eqnarray}
where $(q,r)$ satisfies $2\leq q,r \leq \infty$ and 3/q=1/2-1/r.
\end{lemma}

As was said in the introduction, $F^s$ can be embedded into many
space-time spaces. We prove a variant version of Lemma 4.1 in
\cite{Taoscatteringgkdv}.
\begin{lemma}[Extension lemma]\label{lemmaextension}
Let $Y$ be any space-time Banach space which obeys the time
modulation estimate
\begin{eqnarray}
\norm{g(t)F(t,x)}_{Y}\leq \norm{g}_{L_t^\infty} \norm{F(t,x)}_Y
\end{eqnarray}
for any $F\in Y$ and $g\in L_t^\infty$. Moreover, if for all
$u_{0}\in L_x^2$
\[\norm{W(t)u_{0}}_{Y}\les \norm{u_{0}}_{L_x^2}.\]
Then one also has the estimate that for all $k\in \Z_+$ and $u\in
{F}^0$
\[\norm{P_{k}(u)}_{Y}\les \norm{\wh{P_{k}(u)}}_{{X}_{k}}.\]
\end{lemma}
\begin{proof}
Fix $k\in \Z_+$ and let $u\in F^0$, then we have
\begin{eqnarray*}
P_ku &=&\int \eta_k(\xi)\ft u(\xi,\tau) e^{ix\xi}e^{it\tau}d\xi d\tau \nonumber\\
&=&\sum_{j=0}^\infty\int \eta_j(\tau-\xi^3)\eta_k(\xi)\ft
u(\xi,\tau) e^{ix\xi}e^{it\tau}d\xi d\tau \nonumber\\
&=&\sum_{j=0}^\infty\int \eta_j(\tau)e^{it\tau} \int \eta_k(\xi)\ft
u(\xi,\tau+\xi^3) e^{ix\xi}e^{it\xi^3}d\xi d\tau.
\end{eqnarray*}
Then from the assumption we get
\begin{eqnarray*}
\norm{P_k(u)}_{Y}\les \sum_{j=0}^\infty \int_{\R} \eta_j(\tau)
\norm{\eta_k(\xi)\ft u(\xi,\tau+\xi^3)}_{L^2_\xi}d\tau\les
\norm{\wh{P_{k}(u)}}_{{X}_{k}}.
\end{eqnarray*}
Therefore, we complete the proof of the lemma.
\end{proof}

Using  Lemma \ref{lemmakdvfree} and Lemma \ref{lemmaextension} we
immediately get
\begin{proposition}[$X_k$ embedding]\label{propXkembedding}
Let $k\in \Z_+$, $j\in \N$ and $(q,r)$ as in Lemma
\ref{lemmakdvfree}. Assume $u\in F^0$, then we have
\begin{eqnarray}
&&\norm{P_k(u)}_{L_t^qL_x^r}\les \norm{\ft[P_k(u)]}_{X_k},\\
&&\norm{P_k(u)}_{L_x^2L_{t\in I}^\infty}\les 2^{3k/4}
\norm{\ft[P_k(u)]}_{X_k},\\
&&\norm{P_k(u)}_{L_x^4L_{t}^\infty}\les 2^{k/4}
\norm{\ft[P_k(u)]}_{X_k},\\
&&\norm{P_j(u)}_{L_x^\infty L_{t}^2}\les
2^{-j}\norm{\ft[P_j(u)]}_{X_j},
\end{eqnarray}
As a consequence, we get from the definition that for $u\in
\bar{F}^s$
\[\norm{u}_{L_t^\infty H^s}\les \norm{u}_{\bar{F}^s}.\]
\end{proposition}

For $k\in \Z$ and $j\in \Z_+$ we define
\[D_{k,j}=\{(\xi,\tau): \xi \in [2^{k-1},2^{k+1}] \mbox{ and } \tau-\xi^3\in I_j\}.\]
Following the $[k;Z]$ methods \cite{Taokz} the bilinear estimates in
$X^{s,b}$ space reduce to some dyadic summations and estimates on
the operator norm: for any $k_1,k_2,k_3\in \Z$ and $j_1,j_2,j_3\in
\Z_+$
\begin{eqnarray}\label{eq:3zmult}
\sup_{(u_{k_2,j_2},\ v_{k_3,j_3})\in
E}\norm{1_{D_{k_1,j_1}}(\xi,\tau)\cdot
u_{k_2,j_2}*v_{k_3,j_3}(\xi,\tau)}_{L_{\xi,\tau}^2}
\end{eqnarray}
where the supremum is taken over on $E$
\[E=\{(u,v):\norm{u}_2,\ \norm{v}_2\leq
1 \mbox{ and } \supp(u) \subset D_{k_2,j_2},\ \supp(v) \subset
D_{k_3,j_3}\}.\] By checking the support properties, we get that in
order for \eqref{eq:3zmult} to be nonzero one must have
\begin{eqnarray}
&& |k_{max}-k_{med}|\leq 3,\label{eq:kmaxkmin}\\
&& 2^{j_{max}}\sim \max(2^{j_{med}}, 2^{k_{max}^2
k_{min}}).\label{eq:LmaxH}
\end{eqnarray}
Sharp estimates on \eqref{eq:3zmult} were obtained in \cite{Taokz}.
We will use these to prove the dyadic bilinear estimates.

\begin{proposition}[Proposition 6.1, \cite{Taokz}]\label{pchar}
Let $k_1,k_2,k_3 \in \Z$ and $j_1,j_2,j_3\in \Z_+$ obey
\eqref{eq:kmaxkmin} and \eqref{eq:LmaxH}. Let $N_i=2^{k_i}$ and
$L_i=2^{j_i}$ for $i=1,2,3$. Then

(i) If $N_{max}\sim N_{min}$ and $L_{max}\sim N_{max}^2N_{min}$,
then we have
\begin{equation}\label{eq:chari}
\eqref{eq:3zmult} \les L_{min}^{1/2}N_{max}^{-1/4}L_{med}^{1/4}.
\end{equation}

(ii) If $N_2\sim N_3 \gg N_1$ and $N_{max}^2N_{min}\sim L_1\ges
L_2,L_3$, then
\begin{equation}\label{eq:charii}
\eqref{eq:3zmult} \les
L_{min}^{1/2}N_{max}^{-1}\min(N_{max}^2N_{min},\frac{N_{max}}{N_{min}}L_{med})^{1/2}.
\end{equation}
Similarly for permutations.

(iii) In all other cases, we have
\begin{equation}\label{eq:chariii}
\eqref{eq:3zmult} \les
L_{min}^{1/2}N_{max}^{-1}\min(N_{max}^2N_{min},L_{med})^{1/2}.
\end{equation}
\end{proposition}

Now we are ready to prove our dyadic bilinear estimates. The first
case is high-low interactions

\begin{proposition}[high-low]\label{prophighlow}
(a) If $k\geq 10$, $|k-k_2|\leq 5$, then for any $u,v\in \bar{F}^s$
\begin{eqnarray}\label{eq:highlowa}
\norm{(i+\tau-\xi^3)^{-1}\eta_k(\xi)i\xi\widehat{P_{\leq
0}u}*\widehat{P_{k_2}v}}_{X_k}\les \norm{{P_{\leq
0}u}}_{L_x^2L_t^\infty}\norm{\widehat{P_{k_2}v}}_{X_{k_2}}.
\end{eqnarray}

(b) If $k\geq 10$, $|k-k_2|\leq 5$ and $1\leq k_1\leq k-9$. Then for
any $u, v\in \bar{F}^s$
\begin{eqnarray}\label{eq:highlowb}
\norm{(i+\tau-\xi^3)^{-1}\eta_k(\xi)i\xi\widehat{P_{k_1}u}*\widehat{P_{k_2}v}}_{X_k}\les\
k^32^{-k/2}2^{-k_1}
\norm{\widehat{P_{k_1}u}}_{X_{k_1}}\norm{\widehat{P_{k_2}v}}_{X_{k_2}}.
\end{eqnarray}
\end{proposition}
\begin{proof}
For simplicity of notations we assume $k=k_2$. For part (a), it
follows from the definition of $X_k$ that
\begin{eqnarray}
\norm{(i+\tau-\xi^3)^{-1}\eta_k(\xi)i\xi\widehat{P_0u}*\widehat{P_kv}}_{X_k}\les
2^k\sum_{j\geq
0}2^{-j/2}\norm{\widehat{P_0u}*\widehat{P_{k_2}v}}_{L_{\xi,\tau}^2}.
\end{eqnarray}
From Plancherel's equality and Proposition \ref{propXkembedding} we
get
\[2^k\norm{\widehat{P_0u}*\widehat{P_{k_2}v}}_{L_{\xi,\tau}^2}\les
2^{k}\norm{{P_0u}}_{L_x^2L_t^\infty}\norm{{P_ku}}_{L_x^\infty
L_t^2}\les
\norm{{P_0u}}_{L_x^2L_t^\infty}\norm{\widehat{P_{k}v}}_{X_{k}},\]which
is part (a) as desired. For part (b), from the definition we get
\begin{eqnarray}\label{eq:highlowb1}
\norm{(i+\tau-\xi^3)^{-1}\eta_k(\xi)i\xi\widehat{P_{k_1}u}*\widehat{P_kv}}_{X_k}\les2^k\sum_{j_i\geq
0}2^{-j_3/2}\norm{1_{D_{k,j_3}}\cdot u_{k_1,j_1}*v_{k,j_2}}_2,
\end{eqnarray}
where
\begin{eqnarray}\label{eq:decomuv}
u_{k_1,j_1}=\eta_{k_1}(\xi)\eta_{j_1}(\tau-\xi^3)\widehat{u},\
v_{k,j_2}=\eta_k(\xi)\eta_{j_2}(\tau-\xi^3)\widehat{v}.
\end{eqnarray}
From \eqref{eq:LmaxH} we may assume $j_{max}\geq 2k+k_1-10$ in the
summation on the right-hand side of \eqref{eq:highlowb1}. We may
also assume $j_1,j_2,j_3\leq 10k$, since otherwise we will apply the
trivial estimates
\[\norm{1_{D_{k_3,j_3}}\cdot u_{k_1,j_1}*v_{k,j_2}}_2\les 2^{j_{min}/2}2^{k_{min}/2}\norm{u_{k_1,j_1}}_2\norm{u_{k_2,j_2}}_2,\]
then there is a $2^{-5k}$ to spare which suffices to give the bound
\eqref{eq:highlowb}. Thus by applying \eqref{eq:charii} we get
\begin{eqnarray}
&&2^k\sum_{j_3,j_1,j_2\geq 0}2^{-j_3/2}\norm{1_{D_{k,j}}u_{k_1,j_1}*v_{k,j_2}}_2\nonumber\\
&&\les \ 2^k\sum_{j_3,j_1,j_2\geq 0}2^{-j/2}2^{j_{min}/2}2^{-k/2}2^{-k_1/2}2^{j_{med}/2}\norm{u_{k_1,j_1}}_2\norm{v_{k,j_2}}_2\nonumber\\
&&\les \ 2^k\sum_{j_{max}\geq
2k+k_1-10}k^32^{-k/2}2^{-k_1/2}2^{-j_{max}/2}\norm{\widehat{P_{k_1}u}}_{X_{k_1}}\norm{\widehat{P_kv}}_{X_k}\nonumber\\
&&\les \
k^32^{-k/2}2^{-k_1}\norm{\widehat{P_{k_1}u}}_{X_{k_1}}\norm{\widehat{P_kv}}_{X_k},
\end{eqnarray}
which completes the proof of the proposition.
\end{proof}

In \cite{GW} we proved a similar result as part (a) but with
$\norm{{P_0u}}_{L_x^2L_t^\infty}$ replaced by
$\norm{\widehat{P_{0}u}}_{X_{0}}$ on the right-hand side of
\eqref{eq:highlowa}. Then we see from \eqref{eq:xbar0x0} that the
high-low interactions are still under control if we assume a little
weaker structure on the low frequency. When the low frequency is
comparable to the high frequency, then we have the following

\begin{proposition}\label{prophhh}
If $k\geq 10$, $|k-k_2|\leq 5$ and $k-9\leq k_1\leq k+10$, then for
any $u,\ v \in F^{-3/4}$
\begin{eqnarray}
\norm{(i+\tau-\xi^3)^{-1}\eta_{k_1}(\xi)i\xi\widehat{P_{k}u}*\widehat{P_{k_2}v}}_{X_{k_1}}\les
\  2^{-3k/4}
\norm{\widehat{P_{k}u}}_{X_{k}}\norm{\widehat{P_{k_2}v}}_{X_{k_2}}.
\end{eqnarray}
\end{proposition}

\begin{proof}
As in the proof of Proposition \ref{prophighlow} we assume $k=k_2$
and it follows from the definition of $X_{k_1}$ that
\begin{eqnarray}\label{eq:highhigh2}
&&\norm{(i+\tau-\xi^3)^{-1}\eta_{k_1}(\xi)i\xi\widehat{P_{k}u}*\widehat{P_kv}}_{X_{k_1}}\nonumber\\
&&\les \ 2^{k_1}\sum_{j_1,j_2,j_3\geq
0}2^{-j_1/2}\norm{1_{D_{k_1,j_1}}u_{k,j_2}*v_{k,j_3}}_2,
\end{eqnarray}
where $u_{k,j_1},v_{k,j_2}$ are as in \eqref{eq:decomuv} and we may
assume $j_{max}\geq 3k-20$ and $j_1,j_2,j_3\leq 10k$ in the
summation. Applying \eqref{eq:chari} we get
\begin{eqnarray*}
&&2^{k_1}\sum_{j_1,j_2,j_3\geq
0}2^{-j_1/2}\norm{1_{D_{k_1,j_1}}u_{k,j_2}*v_{k,j_3}}_2\\
&&\les \big(\sum_{j_1=j_{max}}+\sum_{j_2=j_{max}}+\sum_{j_3=j_{max}}\big)2^{-j_1/2}2^{3k/4}2^{j_{min}/2}2^{j_{med}/4}\norm{u_{k,j_2}}_2\norm{v_{k,j_3}}_2\\
&&:=I+II+III.
\end{eqnarray*}
For the contribution of $I$, since it is easy to get the bound, thus
we omit the details. We only need to bound $II$ in view of the
symmetry. We get that
\begin{eqnarray*}
II&\les&(\sum_{j_2=j_{max},j_1\leq j_3}+\sum_{j_2=j_{max},j_1\geq
j_3})2^{-j_1/2}2^{3k/4}2^{j_{min}/2}2^{j_{med}/4}\norm{u_{k,j_2}}_2\norm{v_{k,j_3}}_2\\
&:=&II_1+II_2.
\end{eqnarray*}
For the contribution of $II_1$, by summing on $j_1$ we have
\begin{eqnarray*}
II_1&\les&\sum_{j_2=j_{max},j_1\leq j_3}2^{-j_1/2}2^{3k/4}2^{j_{1}/2}2^{j_{3}/4}\norm{u_{k,j_2}}_2\norm{v_{k,j_3}}_2\\
&\les&\sum_{j_2\geq 3k-20,j_3\geq
0}2^{3k/4}2^{j_{3}/2}\norm{u_{k,j_2}}_2\norm{v_{k,j_3}}_2\\
&\les&2^{-3k/4}
\norm{\widehat{P_{k}u}}_{X_{k}}\norm{\widehat{P_{k_2}v}}_{X_{k_2}},
\end{eqnarray*}
which is acceptable. For the contribution of $II_2$, we have
\begin{eqnarray*}
II_2&\les&\sum_{j_2=j_{max},j_1\geq j_3}2^{-j_1/2}2^{3k/4}2^{j_{3}/2}2^{j_{1}/4}\norm{u_{k,j_2}}_2\norm{v_{k,j_3}}_2\\
&\les&2^{-3k/4}
\norm{\widehat{P_{k}u}}_{X_{k}}\norm{\widehat{P_{k_2}v}}_{X_{k_2}}.
\end{eqnarray*}
Therefore, we complete the proof of the proposition.
\end{proof}

We consider now $low\times low \ra low$ interaction. Generally
speaking, this case is always easy to handle in many situations.
\begin{proposition}[low-low]\label{proplowlow}
If $0\leq k_1,k_2,k_3\leq 100$, then for any $u,\ v \in F^s$
\begin{eqnarray}\label{eq:lowlow}
\norm{(i+\tau-\xi^3)^{-1}\eta_{k_1}(\xi)i\xi\widehat{\psi(t)P_{k_2}(u)}*\widehat{P_{k_3}(v)}}_{X_{k_1}}\les
\norm{{P_{k_2}u}}_{L_t^\infty L_x^2}\norm{{P_{k_3}v}}_{L_t^\infty
L_x^2}.
\end{eqnarray}
\end{proposition}
\begin{proof}
From the definition of $X_{k_1}$, Plancherel's equality and
Bernstein's inequality we get that
\begin{eqnarray*}
&&\norm{(i+\tau-\xi^3)^{-1}\eta_{k_1}(\xi)i\xi\ft[{\psi(t)P_{k_2}u}]*\ft[{P_{k_3}v}](\xi,\tau)}_{X_{k_1}}\\
&&\les \ 2^{k_1}\sum_{j_3\geq 0}2^{-j_3/2}\norm{\psi(t)P_{k_2}u\cdot
P_{k_3}v}_{L_t^2L_x^2}\\
&&\les \norm{{P_{k_2}u}}_{L_t^\infty
L_x^2}\norm{{P_{k_3}v}}_{L_t^\infty L_x^2},
\end{eqnarray*}
which completes the proof of the Proposition.
\end{proof}

The final case is $high\times high \rightarrow low$ interactions. It
is easy to see that this case is the worst, since $s<0$ and
$\norm{u}_{F^s}, \norm{v}_{F^s}$ are small for $u,v$ with very high
frequency.
\begin{proposition}[high-high]\label{prophighhigh}
(a) If $k\geq 10$, $|k-k_2|\leq 5$, then for any $u,\ v \in F^s$
\begin{eqnarray}\label{eq:highhigha}
\norm{(i+\tau-\xi^3)^{-1}\eta_0(\xi)i\xi\widehat{P_{k}u}*\widehat{P_{k_2}v}}_{X_0}\les
\  k2^{-3k/2}
\norm{\widehat{P_{k}u}}_{X_{k}}\norm{\widehat{P_{k_2}v}}_{X_{k_2}}.
\end{eqnarray}

(b) If $k\geq 10$, $|k-k_2|\leq 5$ and $1\leq k_1\leq k-9$, then for
any $u,\ v \in F^s$
\begin{eqnarray*}
\norm{(i+\tau-\xi^3)^{-1}\eta_{k_1}(\xi)i\xi\widehat{P_{k}u}*\widehat{P_{k_2}v}}_{X_{k_1}}\les
(2^{-3k/2}+k2^{-2k+\half{k_1}})
\norm{\widehat{P_{k}u}}_{X_{k}}\norm{\widehat{P_{k_2}v}}_{X_{k_2}}.
\end{eqnarray*}
\end{proposition}
\begin{proof}
For part (a), as before we assume $k=k_2$ and from the definition we
get the left-hand side of \eqref{eq:highhigha} is dominated by
\begin{eqnarray}
\sum_{k_3=-\infty}^0 2^{k_3}\sum_{j_1,j_2,j_3\geq
0}2^{-j_3/2}\norm{1_{D_{k_3,j_3}}\cdot u_{k,j_1}*v_{k,j_2}}_2,
\end{eqnarray}
where $u_{k,j_1}, v_{k,j_2}$ are as in \eqref{eq:decomuv} and we may
assume that $k_3\geq -10k$ and $j_1,j_2,j_3\leq 10k$. It suffices to
consider the worst case $|j_3-2k-k_3|\leq 10$. Then applying
\eqref{eq:charii} we get that
\begin{eqnarray}
&&\norm{(i+\tau-\xi^3)^{-1}\eta_0(\xi)i\xi\widehat{P_{k}u}*\widehat{P_kv}}_{X_0}\nonumber\\
&&\les \ \sum_{k_3=-10k}^0\sum_{j_1,j_2\geq 0}2^{-k}2^{-k_3/2}
2^{k_3}2^{-k/2}2^{-k_3/2}2^{j_{1}/2}2^{j_{2}/2}\norm{u_{k,j_1}}_2\norm{v_{k,j_2}}_2\nonumber\\
&&\les
k2^{-3k/2}\norm{\widehat{P_{k}u}}_{X_{k}}\norm{\widehat{P_kv}}_{X_k},
\end{eqnarray}
which is part (a). For part (b) we assume $k=k_2$ and it follows
from the definition of $X_{k_1}$ that
\begin{eqnarray}\label{eq:highhigh2}
&&\norm{(i+\tau-\xi^3)^{-1}\eta_{k_1}(\xi)i\xi\widehat{P_{k}u}*\widehat{P_kv}}_{X_{k_1}}\nonumber\\
&&\les \ 2^{k_1}\sum_{j_1,j_2,j_3\geq
0}2^{-j_1/2}\norm{1_{D_{k_1,j_1}}u_{k,j_2}*v_{k,j_3}}_2,
\end{eqnarray}
where $u_{k,j_2},v_{k,j_3}$ are as in \eqref{eq:decomuv}. For the
same reasons as in the proof of Proposition \ref{prophighlow} we may
assume $j_{max}\geq 2k+k_1-10$ and $j_1,j_2,j_3\leq 10k$. We will
bound the right-hand side of \eqref{eq:highhigh2} case by case. The
first case is that $j_1=j_{max}$ in the summation. Then we apply
\eqref{eq:charii} and get that
\begin{eqnarray*}
&& 2^{k_1}\sum_{j_1,j_2,j_3\geq
0}2^{-j_1/2}\norm{1_{D_{k_1,j_1}}u_{k,j_2}*v_{k,j_3}}_2\nonumber\\
&&\les \ 2^{k_1}\sum_{j_{1}\geq 2k+k_1-10}\sum_{j_2,j_3\geq
0}2^{-j_1/2}2^{-k/2}2^{{-k_1}/2}2^{(j_{2}+j_3)/2}\norm{u_{k,j_2}}_{2}\norm{v_{k,j_3}}_{2}\\
&&\les
2^{-3k/2}\norm{\widehat{P_{k}u}}_{X_{k}}\norm{\widehat{P_{k_2}v}}_{X_{k_2}},
\end{eqnarray*}
which is acceptable. If $j_2=j_{max}$, then in this case we have
better estimate for the characterization multiplier. By applying
\eqref{eq:chariii} we get
\begin{eqnarray*}
&& 2^{k_1}\sum_{j_1,j_2,j_3\geq
0}2^{-j_1/2}\norm{1_{D_{k_1,j_1}}u_{k,j_2}*v_{k,j_3}}_2\nonumber\\
&&\les \ 2^{k_1}\sum_{j_{2}\geq 2k+k_1-10}\sum_{j_1,j_3\geq
0}2^{-j_1/2}2^{-k}2^{(j_{1}+j_3)/2}\norm{u_{k,j_2}}_{2}\norm{v_{k,j_3}}_{2}\\
&&\les k
2^{-2k}2^{k_1/2}\norm{\widehat{P_{k}u}}_{X_{k}}\norm{\widehat{P_{k_2}v}}_{X_{k_2}},
\end{eqnarray*}
where in the last inequality we use $j_1\leq 10k$. The last case
$j_3=j_{max}$ is identical to the case $j_2=j_{max}$ from symmetry.
Therefore, we complete the proof of the proposition.
\end{proof}

The main reason for us applying $\bar{F}^{-3/4}$ is the logarithmic
loss of derivative in \eqref{eq:highhigha}. We believe that this
loss is essential. Precisely, we conjecture the following: There
doesn't exist a constant $C>0$ such that for all $k\in \N$ and
$u,v\in F^0$
\begin{eqnarray}
\norm{(i+\tau-\xi^3)^{-1}\eta_0(\xi)i\xi\widehat{P_{k}u}*\widehat{P_{k}v}}_{X_0}\leq
C2^{-3k/2}
\norm{\widehat{P_{k}u}}_{X_{k}}\norm{\widehat{P_{k}v}}_{X_{k}}.
\end{eqnarray}
We can't prove it so far. But fortunately we can avoid the
logarithmic loss in \eqref{eq:highhigha} by using a $\bar{X}_0$
structure on the low frequency\footnote{The author is grateful to
Zhaohui Huo for pointing out an error in the proof of this
proposition in an early version of this paper that the general
bilinear estimates do not follow directly from the extension lemma.
Thus we give a direct proof.}.

\begin{proposition}[$\bar{X}_0$ estimate]\label{propXbares}
Let $|k_1-k_2|\leq 5$ and $k_1\geq 10$. Then we have for all $u,v\in
\bar{F}^0$
\begin{eqnarray*}
&\normo{\psi(t)\int_0^t W(t-s)P_{\leq
0}\partial_x[P_{k_1}u(s)P_{k_2}v(s)]ds}_{L_x^2L_t^\infty}\les
2^{-\half{3k_1}}\norm{\wh{P_{k_1}u}}_{X_{k_1}}\norm{\wh{P_{k_2}u}}_{X_{k_2}}.&
\end{eqnarray*}
\end{proposition}

\begin{proof}
Denote $Q(u,v)=\psi(t)\int_0^t W(t-s)P_{\leq
0}\partial_x[P_{k_1}u(s)P_{k_2}v(s)]ds$. By straightforward
computations we get
\begin{eqnarray*}
\ft\left[Q(u,v)\right](\xi,\tau)&=&c\int_\R
\frac{\widehat{\psi}(\tau-\tau')-\widehat{\psi}(\tau-\xi^3)}{\tau'-\xi^3}\eta_0(\xi)i\xi
\\
&&\times\
d\tau'\int_{\xi=\xi_1+\xi_2,\tau'=\tau_1+\tau_2}\wh{P_{k_1}u}(\xi_1,\tau_1)\wh{P_{k_2}v}(\xi_2,\tau_2).
\end{eqnarray*}
Fixing $\xi \in \R$, we decomposing the hyperplane as following
\begin{eqnarray*}
A_1&=&\{\xi=\xi_1+\xi_2,\tau'=\tau_1+\tau_2: |\xi|\les 2^{-2k_1}\};\\
A_2&=&\{\xi=\xi_1+\xi_2,\tau'=\tau_1+\tau_2: |\xi|\gg 2^{-2k_1},
|\tau_i-\xi_i^3|\ll 3\cdot 2^{2k_1}|\xi|, i=1,2\};\\
A_3&=&\{\xi=\xi_1+\xi_2,\tau'=\tau_1+\tau_2: |\xi|\gg 2^{-2k_1},
|\tau_1-\xi_1^3|\ges 3\cdot 2^{2k_1}|\xi|\};\\
A_4&=&\{\xi=\xi_1+\xi_2,\tau'=\tau_1+\tau_2: |\xi|\gg 2^{-2k_1},
|\tau_2-\xi_2^3|\ges 3\cdot 2^{2k_1}|\xi|\}.
\end{eqnarray*}
Then we get
\[\ft\left[\psi(t)\cdot
\int_0^tW(t-s)P_{\leq
0}\partial_x[P_{k_1}u(s)P_{k_2}v(s)]ds\right](\xi,\tau)=I+II+III,\]
where
\begin{eqnarray*}
I&=&C\int_\R
\frac{\widehat{\psi}(\tau-\tau')-\widehat{\psi}(\tau-\xi^3)}{\tau'-\xi^3}\eta_0(\xi)i\xi
\int_{A_1}\wh{P_{k_1}u}(\xi_1,\tau_1)\wh{P_{k_2}v}(\xi_2,\tau_2)d\tau',\\
II&=&C\int_\R
\frac{\widehat{\psi}(\tau-\tau')-\widehat{\psi}(\tau-\xi^3)}{\tau'-\xi^3}\eta_0(\xi)i\xi
\int_{A_2}\wh{P_{k_1}u}(\xi_1,\tau_1)\wh{P_{k_2}v}(\xi_2,\tau_2)d\tau',\\
II&=&C\int_\R
\frac{\widehat{\psi}(\tau-\tau')-\widehat{\psi}(\tau-\xi^3)}{\tau'-\xi^3}\eta_0(\xi)i\xi
\int_{A_3\cup
A_4}\wh{P_{k_1}u}(\xi_1,\tau_1)\wh{P_{k_2}v}(\xi_2,\tau_2)d\tau'.
\end{eqnarray*}

We consider first the contribution of the term $I$. Using
Proposition \ref{propXkembedding} and Proposition \ref{proplineares}
(b), we get
\[\norm{\ft^{-1}(I)}_{L_x^2L_t^\infty}\les \norm{I}_{X_0}\les \normo{(i+\tau'-\xi^3)^{-1}\eta_0(\xi)i\xi
\int_{A_1}\wh{P_{k_1}u}(\xi_1,\tau_1)\wh{P_{k_2}v}(\xi_2,\tau_2)}_{X_0}.\]
Since in the area $A_1$ we have $|\xi|\les 2^{-2k_1}$, thus we get
\begin{eqnarray*}
&&\normo{(i+\tau'-\xi^3)^{-1}\eta_0(\xi)i\xi
\int_{A_1}\wh{P_{k_1}u}(\xi_1,\tau_1)\wh{P_{k_2}v}(\xi_2,\tau_2)}_{X_0}\\
&\les& \sum_{k_3\leq -2k_1+10}\sum_{j_3\geq 0}
2^{-j_3/2}2^{k_3}\sum_{j_1\geq 0, j_2\geq
0}\norm{1_{D_{k_3,j_3}}\cdot u_{k_1,j_1}*v_{k_2,j_2}}_{L^2}
\end{eqnarray*}
where
\begin{eqnarray}
u_{k_1,j_1}(\xi,\tau)=\eta_{k_1}(\xi)\eta_{j_1}(\tau-\xi^3)\wh{u}(\xi,\tau),
v_{k_1,j_1}(\xi,\tau)=\eta_{k_1}(\xi)\eta_{j_1}(\tau-\xi^3)\wh{v}(\xi,\tau).
\end{eqnarray}
Using Proposition \ref{pchar} (iii), then we get
\begin{eqnarray*}
\norm{\ft^{-1}(I)}_{L_x^2L_t^\infty}&\les& \sum_{k_3\leq
-2k_1+10}\sum_{j_i\geq 0} 2^{-j_3/2}2^{k_3}2^{j_{min}/2}2^{k_3/2}
\norm{u_{k_1,j_1}}_{L^2}\norm{v_{k_2,j_2}}_{L^2}\\
&\les&
2^{-3k_1}\norm{\wh{P_{k_1}u}}_{X_{k_1}}\norm{\wh{P_{k_2}u}}_{X_{k_2}},
\end{eqnarray*}
which suffices to give the bound for the term $I$.

Next we consider the contribution of the term $III$. As for the term
$I$, Using Proposition \ref{propXkembedding} and Proposition
\ref{proplineares} (b), we get
\begin{eqnarray*}
\norm{\ft^{-1}(III)}_{L_x^2L_t^\infty}&\les&
\normo{(i+\tau'-\xi^3)^{-1}\eta_0(\xi)i\xi \int_{A_3\cup
A_4}\wh{P_{k_1}u}(\xi_1,\tau_1)\wh{P_{k_2}v}(\xi_2,\tau_2)}_{X_0}\\
&\les&\sum_{k_3\leq 0}\sum_{j_3\geq 0}
2^{-j_3/2}2^{k_3}\sum_{j_1\geq 0, j_2\geq
0}\norm{1_{D_{k_3,j_3}}\cdot u_{k_1,j_1}*v_{k_2,j_2}}_{L^2}
\end{eqnarray*}
Clearly we may assume $j_3\leq 10k_1$ in the summation above.
Without loss of generality, we assume $|\tau_1-\xi_1^3|\ges
3|\xi\xi_1\xi_2|$. Using Proposition \ref{pchar} (iii), then we get
\begin{eqnarray*}
\norm{\ft^{-1}(III)}_{L_x^2L_t^\infty}&\les& \sum_{k_3\leq
0}\sum_{j_1\geq k_3+2k_1-10, j_2,j_3\geq 0}
2^{k_3}2^{j_{2}/2}2^{-k_1}
\norm{u_{k_1,j_1}}_{L^2}\norm{v_{k_2,j_2}}_{L^2}\\
&\les&
k_12^{-2k_1}\norm{\wh{P_{k_1}u}}_{X_{k_1}}\norm{\wh{P_{k_2}u}}_{X_{k_2}},
\end{eqnarray*}
which suffices to give the bound for the term $III$.

Now we consider the contribution of the term $II$. From the proof of
the dyadic bilinear estimates, we know this term is the main
contribution. By computation we get
\begin{eqnarray*}
\ft_t^{-1}(II)&=&\psi(t)\int_0^t e^{i(t-s)\xi^3}\eta_0(\xi)i\xi
\int_{\R^2}e^{is(\tau_1+\tau_2)}\\
&&\times \
\int_{\xi=\xi_1+\xi_2}{u_{k_1}}(\xi_1,\tau_1){v_{k_2}}(\xi_2,\tau_2)\
d\tau_1 d\tau_2ds
\end{eqnarray*}
where
\begin{eqnarray*}
u_{k_1}(\xi_1,\tau_1)&=&\eta_{k_1}(\xi_1)1_{\{|\tau_1-\xi_1^3|\ll
3\cdot 2^{2k_1}|\xi|\}}\wh{u}(\xi_1,\tau_1),\\
v_{k_2}(\xi_2,\tau_2)&=&\eta_{k_2}(\xi_2)1_{\{|\tau_2-\xi_2^3|\ll
3\cdot 2^{2k_1}|\xi|\}}\wh{v}(\xi_2,\tau_2).
\end{eqnarray*}
By a change of variable $\tau_1'=\tau_1-\xi_1^3$,
$\tau_2'=\tau_2-\xi_2^3$, we get
\begin{eqnarray*}
\ft_t^{-1}(II)&=&\psi(t)e^{it\xi^3}\eta_0(\xi)i\xi\int_0^t
e^{-is\xi^3}
\int_{\R^2}e^{is(\tau_1+\tau_2)}\\
&&\times \
\int_{\xi=\xi_1+\xi_2}e^{is\xi_1^3}{u_{k_1}}(\xi_1,\tau_1+\xi_1^3)e^{is\xi_2^3}{v_{k_2}}(\xi_2,\tau_2+\xi_2^3)\
d\tau_1 d\tau_2ds\\
&=&\psi(t)e^{it\xi^3}\eta_0(\xi)\xi\int_{\R^2}e^{it(\tau_1+\tau_2)}\int_{\xi=\xi_1+\xi_2}\frac{e^{it(\xi_1^3+\xi_2^3-\xi^3)}-e^{-it(\tau_1+\tau_2)}}{\tau_1+\tau_2-\xi^3+\xi_1^3+\xi_2^3}\\
&&\times \
{u_{k_1}}(\xi_1,\tau_1+\xi_1^3){v_{k_2}}(\xi_2,\tau_2+\xi_2^3)\
d\tau_1 d\tau_2\\
&:=&\ft_t^{-1}(II_1)-\ft_t^{-1}(II_2).
\end{eqnarray*}
For the contribution of the term $II_2$, we have
\begin{eqnarray*}
\ft_t^{-1}(II_2)=\int_{\R^2}\psi(t)e^{it\xi^3}\eta_0(\xi)\xi\int_{\xi=\xi_1+\xi_2}
\frac{{u_{k_1}}(\xi_1,\tau_1+\xi_1^3){v_{k_2}}(\xi_2,\tau_2+\xi_2^3)}{\tau_1+\tau_2-\xi^3+\xi_1^3+\xi_2^3}\
d\tau_1 d\tau_2.
\end{eqnarray*}
Since in the integral area we have
$|\tau_1+\tau_2-\xi^3+\xi_1^3+\xi_2^3|\sim |\xi\xi_1\xi_2|$, thus we
get from Lemma \ref{lemmakdvfree} that
\begin{eqnarray*}
\norm{\ft^{-1}(II_2)}_{L_x^2L_t^\infty}&\les&
\int_{\R^2}\normo{\int_{\xi=\xi_1+\xi_2}\xi\frac{{u_{k_1}}(\xi_1,\tau_1+\xi_1^3){v_{k_2}}(\xi_2,\tau_2+\xi_2^3)}{\tau_1+\tau_2-\xi^3+\xi_1^3+\xi_2^3}}_{L_\xi^2}d\tau_1d\tau_2\\
&\les&2^{-\half{3k_1}}\norm{\wh{P_{k_1}u}}_{X_{k_1}}\norm{\wh{P_{k_2}u}}_{X_{k_2}}.
\end{eqnarray*}

To prove the proposition, it remains to prove the following
\begin{eqnarray*}
\norm{\ft^{-1}(II_1)}_{L_x^2L_t^\infty}\les
2^{-3k_1/2}\norm{\wh{P_{k_1}u}}_{X_{k_1}}\norm{\wh{P_{k_2}u}}_{X_{k_2}}.
\end{eqnarray*}
Compare the term $II_1$ with the following term $II'_1$:
\begin{eqnarray*}
\ft_t^{-1}(II'_1)
&=&\psi(t)e^{it\xi^3}\eta_0(\xi)\xi\int_{\R^2}e^{it(\tau_1+\tau_2)}\int_{\xi=\xi_1+\xi_2}\frac{e^{it(\xi_1^3+\xi_2^3-\xi^3)}}{-\xi^3+\xi_1^3+\xi_2^3}\\
&&\times \
{u_{k_1}}(\xi_1,\tau_1+\xi_1^3){v_{k_2}}(\xi_2,\tau_2+\xi_2^3)\
d\tau_1 d\tau_2.
\end{eqnarray*}
For the term $II'_1$ we have
\begin{eqnarray*}
\ft_t^{-1}(II'_1)
&=&\int_{\R^2}\psi(t)\eta_0(\xi) e^{it(\tau_1+\tau_2)}1_{\{|\xi|\gg |\tau_1|2^{-2k_1}\}}1_{\{|\xi|\gg |\tau_2|2^{-2k_1}\}}\\
&&\times \
\int_{\xi=\xi_1+\xi_2}\frac{e^{it(\xi_1^3+\xi_2^3)}}{-3\xi_1\xi_2}\ft(f_{\tau_1})(\xi_1)\ft(g_{\tau_2})(\xi_2)\
d\tau_1 d\tau_2.
\end{eqnarray*}
where for $\tau_1,\tau_2\in \R$, we set
\[\ft(f_{\tau_1})(\xi)=\wh{P_{k_1}u}(\xi,\tau_1+\xi^3),\
\ft(g_{\tau_2})(\xi)=\wh{P_{k_2}v}(\xi,\tau_2+\xi^3).\] Since it is
easy to see that (actually we need a smooth version of
$1_{\{|\xi|\gg \lambda\}}$): $\forall \ \lambda>0$ we have
\[\norm{\ft_x^{-1}1_{\{|\xi|\gg \lambda\}}\ft_x u}_{L_x^2L_t^\infty}\les \norm{u}_{L_x^2L_t^\infty},\]
thus we get from Lemma \ref{lemmakdvfree} that
\begin{eqnarray*}
\norm{\ft^{-1}(II'_1)}_{L_x^2L_t^\infty}&\les&
\int_{\R^2}\norm{W(t)\partial_x^{-1}f_{\tau_1}W(t)\partial_x^{-1}f_{\tau_2}}_{L_x^2L_t^\infty}d\tau_1d\tau_2\\
&\les&
\int_{\R^2}\norm{W(t)\partial_x^{-1}f_{\tau_1})}_{L_x^4L_t^\infty}\norm{W(t)\partial_x^{-1}f_{\tau_2}}_{L_x^4L_t^\infty}d\tau_1d\tau_2\\
&\les&2^{-\half{3k_1}}\norm{\wh{P_{k_1}u}}_{X_{k_1}}\norm{\wh{P_{k_2}u}}_{X_{k_2}},
\end{eqnarray*}
which gives the bound for the term $II'_1$.

To prove the proposition, it remains to prove the following
\begin{eqnarray*}
\norm{\ft^{-1}(II_1-II'_1)}_{L_x^2L_t^\infty}\les
2^{-3k_1/2}\norm{\wh{P_{k_1}u}}_{X_{k_1}}\norm{\wh{P_{k_2}u}}_{X_{k_2}}.
\end{eqnarray*}
Since in the integral area we have $|\tau_i|\ll 2^{2k_1}|\xi|$,
$i=1,2$, thus we get in the hyperplane
\[\frac{1}{\tau_1+\tau_2-\xi^3+\xi_1^3+\xi_2^3}=\sum_{n=0}^\infty \rev{3\xi\xi_1\xi_2}\brk{\frac{\tau_1+\tau_2}{3\xi\xi_1\xi_2}}^n \]
and then
\[\frac{1}{\tau_1+\tau_2-\xi^3+\xi_1^3+\xi_2^3}-\frac{1}{-\xi^3+\xi_1^3+\xi_2^3}=\sum_{n=1}^\infty \rev{3\xi\xi_1\xi_2}\brk{\frac{\tau_1+\tau_2}{3\xi\xi_1\xi_2}}^n. \]
Thus we get
\begin{eqnarray*}
&&\ft_t^{-1}(II_1-II'_1)\\
&=&\psi(t)\eta_0(\xi)\xi\int_{\R^2}e^{it(\tau_1+\tau_2)}\int_{\xi=\xi_1+\xi_2}\sum_{n=1}^\infty
\rev{3\xi\xi_1\xi_2}\brk{\frac{\tau_1+\tau_2}{3\xi\xi_1\xi_2}}^n\\
&&\times \ e^{it(\xi_1^3+\xi_2^3)}
{u_{k_1}}(\xi_1,\tau_1+\xi_1^3){v_{k_2}}(\xi_2,\tau_2+\xi_2^3)\
d\tau_1 d\tau_2.
\end{eqnarray*}
Then decomposing the low frequency, we get
\begin{eqnarray*}
&&\ft_t^{-1}(II_1-II'_1)\\
&=&\sum_{n=1}^\infty \int_{\R^2}
e^{it(\tau_1+\tau_2)}\sum_{2^{k_3}\gg
2^{-2k_1}\max(|\tau_1|,|\tau_2|)}\psi(t)e^{it(\xi_1^3+\xi_2^3)}\chi_{k_3}(\xi)\\
&&\times \
\int_{\xi=\xi_1+\xi_2}\brk{\frac{\tau_1+\tau_2}{3\xi\xi_1\xi_2}}^n
\frac{{u_{k_1}}(\xi_1,\tau_1+\xi_1^3){v_{k_2}}(\xi_2,\tau_2+\xi_2^3)}{3\xi_1\xi_2}
\ d\tau_1 d\tau_2
\end{eqnarray*}
We rewrite it as
\begin{eqnarray*}
&&\ft_t^{-1}(II_1-II'_1)\\
&=&\sum_{n=1}^\infty \int_{\R^2}
e^{it(\tau_1+\tau_2)}\sum_{2^{k_3}\gg
2^{-2k_1}\max(|\tau_1|,|\tau_2|)}\psi(t)e^{it(\xi_1^3+\xi_2^3)}\chi_{k_3}(\xi)(\xi/2^{k_3})^{-n} \\
&&\times \ 2^{-nk_3}
\int_{\xi=\xi_1+\xi_2}\brk{\frac{\tau_1+\tau_2}{3\xi_1\xi_2}}^n
\frac{{u_{k_1}}(\xi_1,\tau_1+\xi_1^3){v_{k_2}}(\xi_2,\tau_2+\xi_2^3)}{3\xi_1\xi_2}
\ d\tau_1 d\tau_2
\end{eqnarray*}
Using the fact that $\chi_{k_3}(\xi)(\xi/2^{k_3})^{-n}$ is a
multiplier for the space $L_x^2L_t^\infty$ and as for the term
$II_1'$, we get
\begin{eqnarray*}
&&\norm{\ft^{-1}(II_1-II'_1)}_{L_x^2L_t^\infty}\\
&\les &\sum_{n=1}^\infty \int_{\R^2} \sum_{2^{k_3}\gg
2^{-2k_1}\max(|\tau_1|,|\tau_2|)} C^n |\tau_1+\tau_2|^n2^{-nk_3}\\
&&\times \  \normo{\psi(t)e^{it(\xi_1^3+\xi_2^3)}
\int_{\xi=\xi_1+\xi_2}
\frac{\ft(f_{\tau_1})(\xi_1)\ft(g_{\tau_2})(\xi_2)}{3\xi_1^{n+1}\xi_2^{n+1}}}_{L_x^2L_t^\infty}
\ d\tau_1 d\tau_2.
\end{eqnarray*}
Using Lemma \ref{lemmakdvfree} and summing on $k_3$, we get that for
some $M\gg 1$
\begin{eqnarray*}
&&\norm{\ft^{-1}(II_1-II'_1)}_{L_x^2L_t^\infty}\\
&\les &\sum_{n=1}^\infty \int_{\R^2} \sum_{2^{k_3}\gg
2^{-2k_1}\max(|\tau_1|,|\tau_2|)} C^n
|\tau_1+\tau_2|^n2^{-nk_3}2^{-2nk_1}\\
&&\times \
2^{-3k_1/2}\norm{\ft(f_{\tau_1})}_{L^2}\norm{\ft(g_{\tau_2})}_{L^2}d\tau_1
d\tau_2\\
&\les& \sum_{n=1}^\infty (C/M)^n \int_{\R^2}
2^{-3k_1/2}\norm{\ft(f_{\tau_1})}_{L^2}\norm{\ft(g_{\tau_2})}_{L^2}d\tau_1
d\tau_2\\
&\les&
2^{-3k_1/2}\norm{\wh{P_{k_1}u}}_{X_{k_1}}\norm{\wh{P_{k_2}u}}_{X_{k_2}}.
\end{eqnarray*}
Therefore, we complete the proof of the proposition.
\end{proof}

\section{Proof of Theorem \ref{thmlwp} and \ref{thmgwp}}

We construct first a strong local solution to the KdV equation
\eqref{eq:kdv} via contraction principle. The main ingredients are
the dyadic bilinear estimates obtained in the last section. We
observe first that the KdV equation \eqref{eq:kdv} is invariant
under the following scaling transform: for $\lambda>0$
\begin{equation}\label{eq:scaling}
u(x,t)\rightarrow \lambda^2 u(\lambda x, \lambda^3 t),\
\phi(x)\rightarrow \lambda^2\phi(\lambda x).
\end{equation}
$\dot{H}^{-3/2}$ is the critical space to \eqref{eq:kdv} in the
sense that $\norm{\lambda^2\phi(\lambda
\cdot)}_{\dot{H}^{-3/2}}=\norm{\phi}_{\dot{H}^{-3/2}}$. From the
fact that
\[\norm{\lambda^2\phi(\lambda x)}_{H^{-3/4}}\les \lambda^{3/2}\norm{\phi}_{H^{-3/4}}+\lambda^{3/4}\norm{\phi}_{H^{-3/4}}\]
then by taking $\lambda$ sufficiently small we may assume
\begin{eqnarray}\label{eq:smalldata}
\norm{\phi}_{H^{-3/4}}\leq \epsilon \ll 1.
\end{eqnarray}
Then we only need to construct the solution of \eqref{eq:kdv} on
$[-1,1]$ under the condition \eqref{eq:smalldata}. From Duhamel's
principle, \eqref{eq:kdv} is equivalent to the integral equation
\begin{equation}\label{eq:kdvint}
u(t)=W(t)\phi-\half{1}\int_0^tW(t-\tau)\partial_x(u^2(\tau))d\tau.
\end{equation}
We will apply a fixed point argument to solve the following
truncated version
\begin{equation}\label{eq:kdvinttrun}
u(t)=\psi(\frac{t}{4})\left[W(t)\phi-\int_0^tW(t-\tau)\partial_x(\psi^2(\tau)u^2(\tau))d\tau
\right].
\end{equation}
It is easy to see that if $u$ solves \eqref{eq:kdvinttrun} then $u$
is a solution of \eqref{eq:kdvint} and hence of \eqref{eq:kdv} on
the time interval $[-1, 1]$.

\begin{proposition}[Linear estimates]\label{proplineares}
(a) Assume $s\in \R$ and  $\phi \in H^s$. Then there exists $C>0$
such that
\begin{eqnarray}
\norm{\psi(t)W(t)\phi}_{\bar{F}^s}\leq C\norm{\phi}_{H^{s}}.
\end{eqnarray}

(b) Assume $s\in \R, k\in \Z_+$ and $u$ satisfies
$(i+\tau-\xi^3)^{-1}\ft(u)\in X_k$. Then there exists $C>0$ such
that
\begin{eqnarray}
\normo{\ft\left[\psi(t)\int_0^t W(t-s)(u(s))ds\right]}_{X_k}\leq
C\norm{(i+\tau-\xi^3)^{-1}\ft(u)}_{X_k}.
\end{eqnarray}
\end{proposition}

\begin{proof}
Part (a) follows from Proposition \ref{lemmakdvfree} and the
definitions. Part (b) has appeared in many literatures, see for
example \cite{In-Ke,GW}.
\end{proof}

For $u,v\in \bar{F}^s$ we define the bilinear operator
\begin{eqnarray}
B(u,v)=\psi(\frac{t}{4})\int_0^tW(t-\tau)\partial_x\big(\psi^2(\tau)u(\tau)\cdot
v(\tau)\big)d\tau.
\end{eqnarray}
In order to apply a fixed point argument to \eqref{eq:kdvinttrun},
all the issues are then reduced to show the boundness of
$B:\bar{F}^s\times \bar{F}^s\rightarrow \bar{F}^s$.

\begin{proposition}[Bilinear estimates]\label{propbilinearbd}
Assume $-3/4\leq s\leq 0$. Then there exists $C>0$ such that
\begin{eqnarray}\label{eq:bilinearbd}
\norm{B(u,v)}_{\bar{F}^s}\leq
C(\norm{u}_{\bar{F}^s}\norm{v}_{\bar{F}^{-3/4}}+\norm{u}_{\bar{F}^{-3/4}}\norm{v}_{\bar{F}^s})
\end{eqnarray}
hold for any $u,v\in \bar{F}^s$.
\end{proposition}

\begin{proof}
In view of definition, we get
\begin{eqnarray}\label{eq:bilinearFs}
\norm{B(u,v)}_{\bar{F}^s}^2=\norm{P_{\leq
0}B(u,v)}_{\bar{X}_0}^2+\sum_{k_1\geq
1}2^{2k_1s}\norm{\eta_{k_1}(\xi)\ft[B(u,v)]}_{X_{k_1}}^2.
\end{eqnarray}
We consider first the contribution of the second term on the
right-hand side of \eqref{eq:bilinearFs}. By decomposing $u,v$ we
have
\begin{eqnarray}\label{eq:bilinearFs1}
\norm{\eta_{k_1}(\xi)\ft[B(u,v)]}_{X_{k_1}}\les \sum_{k_2,k_3\geq
0}\norm{\eta_{k_1}(\xi)\ft[B(P_{k_2}(u),P_{k_3}(v))]}_{X_{k_1}}.
\end{eqnarray}
From Proposition \ref{proplineares} (b) the right-hand side of
\eqref{eq:bilinearFs1} is dominated by
\begin{eqnarray}\label{eq:bilinearFs11}
\sum_{k_2,k_3\geq
0}\norm{(i+\tau-\xi^3)^{-1}\eta_{k_1}(\xi)i\xi\wh{\psi(t)P_{k_2}u}*\wh{\psi(t)P_{k_3}v})}_{X_{k_1}}.
\end{eqnarray}
From symmetry we assume $k_2\leq k_3$ in \eqref{eq:bilinearFs11}. It
suffices to prove
\begin{eqnarray}\label{eq:bilinearFs12}
&&\bigg(\sum_{k_1\geq 1}2^{2k_1s}\big[\sum_{k_2,k_3\geq
0}\norm{(i+\tau-\xi^3)^{-1}\eta_{k_1}(\xi)i\xi\wh{\psi(t)P_{k_2}u}*\wh{\psi(t)P_{k_3}v})}_{X_{k_1}}\big]^2\bigg)^{1/2}\nonumber\\
&&\ \les \norm{u}_{\bar{F}^{-3/4}}\norm{v}_{\bar{F}^s}.
\end{eqnarray}
If $k_{max}\leq 20$ then applying Proposition \ref{proplowlow} and
from \eqref{eq:propertyXk2} we get that \eqref{eq:bilinearFs11} is
dominated by
\begin{eqnarray}
\sum_{k_{max}\leq 20}\norm{{P_{k_2}u}}_{L_t^\infty
L_x^2}\norm{{P_{k_3}v}}_{L_t^\infty L_x^2},
\end{eqnarray}
which suffices to give the bound \eqref{eq:bilinearFs12} in this
case since it's easy to see that we have $\norm{P_{k}u}_{L_t^\infty
L_x^2}\les \norm{P_{k}u}_{X_k}$ for $k\geq 1$ and
$\norm{P_{k}u}_{L_t^\infty L_x^2}\les \norm{P_{k}u}_{\bar{X}_k}$ for
$k=0$. Assuming $k_{max}\geq 20$ in \eqref{eq:bilinearFs11}, we have
three cases. If $|k_1-k_3|\leq 5, k_2\leq k_1-10$, then applying
Proposition \ref{prophighlow} (a) for $k_2=0$ and (b) for $k_2\geq
1$; If $|k_1-k_3|\leq 5, k_1-9\leq k_2 \leq k_3$, then applying
Proposition \ref{prophhh}; If $|k_2-k_3|\leq 5, 1\leq k_1\leq
k_2-5$, then applying Proposition \ref{prophighhigh} (b). We easily
get the bound \eqref{eq:bilinearFs12} as desired.

To prove Proposition \ref{propbilinearbd}, it remains to prove that
\begin{eqnarray}\label{eq:bilinearFs0}
\norm{B(u,v)}_{\bar{X}_0}\leq
C(\norm{u}_{\bar{F}^s}\norm{v}_{\bar{F}^{-3/4}}+\norm{u}_{\bar{F}^{-3/4}}\norm{v}_{\bar{F}^s}).
\end{eqnarray}
By decomposing $u, v$ as before we obtain
\begin{eqnarray}
\norm{B(u,v)}_{\bar{X}_0}\leq \sum_{k_2,k_3\geq
0}\norm{B(P_{k_2}u,P_{k_3}v)}_{\bar{X}_0}.
\end{eqnarray}
If $\max(k_2,k_3)\leq 10$, then from \eqref{eq:xbar0x0} and
Proposition \ref{proplineares} and Proposition \ref{proplowlow} we
obtain that
\[\norm{B(P_{k_2}u,P_{k_3}v)}_{\bar{X}_0}\les \norm{{P_{k_2}u}}_{L_t^\infty
L_x^2}\norm{{P_{k_3}v}}_{L_t^\infty L_x^2},\] which suffices to give
the bound \eqref{eq:bilinearFs0} in this case. If $\max(k_2,k_3)\geq
10$, then we must have $|k_2-k_3|\leq 5$. Then from Proposition
\ref{propXbares} we have
\begin{eqnarray}
\norm{B(u,v)}_{\bar{X}_0}&\leq& \sum_{|k_2-k_3|\leq 5,\ k_2,k_3\geq
10}2^{-3k_2/2}\norm{\ft(P_{k_2}u)}_{{X}_{k_2}}\norm{\ft(P_{k_3}v)}_{{X}_{k_3}}\nonumber\\
&\les& \norm{u}_{\bar{F}^{-3/4}}\norm{v}_{\bar{F}^{-3/4}}
\end{eqnarray}
which gives \eqref{eq:bilinearFs0} as desired. Thus we complete the
proof of the proposition.
\end{proof}

With a standard argument (see for example, Lemma 4, \cite{Cann}), we
get that there is a unique solution $u$ to \eqref{eq:kdvinttrun}
such that $\norm{u}_{\bar{F}^{-3/4}}\leq C\epsilon_0$. So far, we
have proved Theorem \ref{thmlwp} (a). The rest of Theorem
\ref{thmlwp} also follow from standard argument.

In the rest of this section we prove Theorem \ref{thmgwp}. The
standard way to extend a local solution to a global one is to make
use of the conservation laws. It is well-known that the KdV equation
is completely integrable and hence has infinite conservation laws.
However, there is no conservation laws below $L^2$, and thus one can
not automatically get global well-posedness below $L^2$. J.
Colliander, M. Keel, G. Staffilani, H. Takaoka and T. Tao developed
the modified energy (I-method) to prove global well-posedness below
the energy norm. We adapt I-method to extend the local solution in
$H^{-3/4}$ and refer the readers to \cite{I-method} for many
details. We define $I-operator$ by
\[\wh{If}(\xi)=m(\xi)\wh{f}(\xi),\]
where the multiplier $m(\xi)$ is smooth, monotone, and of the form
for $N\geq 1$
\begin{eqnarray}\label{eq:m}
m(\xi)=\left\{
\begin{array}{r}
1, \quad \quad |\xi|<N,\\
N^{-s}|\xi|^s,\quad  |\xi|>2N.
\end{array}
\right.
\end{eqnarray}
We state a variant local well-posedness result which follows from
slight argument in the last section and from the same reasons as in
\cite{I-method}. This is used to iterate the solution in the
I-method.
\begin{proposition}
Let $-3/4\leq s\leq 0$. Assume $\phi$ satisfies
$\norm{I\phi}_{L^2(\R)}\leq \epsilon_0\ll 1$.  Then there exists a
unique solution $u$ to \eqref{eq:kdv} on $[-1,1]$ such that
\begin{equation}
\norm{Iu}_{\bar{F}^s(1)}\leq C\epsilon_0.
\end{equation}
\end{proposition}

Then it suffices to control $\norm{Iu}_{L^2}$ for all $t$. Let $g:
\R^k \rightarrow \C$ be a function. We say $g$ is symmetric if
$g(\xi_1,\ldots, \xi_k)=g(\sigma(\xi_1,\ldots, \xi_k))$ for all
$\sigma \in S_k$, the group of all permutations on $k$ objects. The
symmetrization of $g$ is the function
\begin{eqnarray}
[g]_{sym}(\xi_1,\xi_2,\ldots, \xi_k)=\rev{k!}\sum_{\sigma\in
S_k}g(\sigma(\xi_1,\xi_2,\ldots,\xi_k)).
\end{eqnarray}
We define a $k-linear$ functional associated to the multiplier $g$
acting on $k$ functions $u_1,\ldots,u_k$,
\begin{eqnarray}
\Lambda_k(g;u_1,\ldots,u_k)=\int_{\xi_1+\ldots+\xi_k=0}g(\xi_1,\ldots,\xi_k)\widehat{u_1}(\xi_1)\ldots
\widehat{u_k}(\xi_k).
\end{eqnarray}
We will often apply $\Lambda_k$ to $k$ copies of the same function
$u$. $\Lambda_k(g;u,\ldots,u)$ may simply be written $\Lambda_k(g)$.
By the symmetry of the measure on hyperplane, we have
$\Lambda_k(g)=\Lambda_k([g]_{sym})$. For $k\in \N$ denote
\[\alpha_k=\xi_1^3+\ldots+\xi_k^3.\] We define the modified energy
$E_I^2(t)$ by
\begin{equation}
E_I^2(t)=\norm{Iu(t)}_{L^2}^2=\Lambda_2(m(\xi_1)m(\xi_2)).
\end{equation}
Form the new modified energy
\begin{eqnarray*}
&&E_I^3(t)=E_I^2(t)+\Lambda_3(\sigma_3),\\
&&E_I^4(t)=E_I^3(t)+\Lambda_4(\sigma_4),
\end{eqnarray*}
where
\begin{eqnarray*}
&&\sigma_3=-\frac{M_3}{\alpha_3},\quad
M_3(\xi_1,\xi_2,\xi_3)=-i[m(\xi_1)m(\xi_2+\xi_3)(\xi_2+\xi_3)]_{sym};\\
&&\sigma_4=-\frac{M_4}{\alpha_4},\quad
M_4(\xi_1,\xi_2,\xi_3,\xi_4)=-i\half
3[\sigma_3(\xi_1,\xi_2,\xi_3+\xi_4)(\xi_3+\xi_4)]_{sym}.
\end{eqnarray*}

\begin{proposition}
Let $I$ be defined with the multiplier $m$ of the form \eqref{eq:m}
and $s=-3/4$. Then
\begin{equation}
|E_I^4(t)-E_I^2(t)|\les \norm{Iu(t)}_{L^2}^3+\norm{Iu(t)}_{L^2}^4.
\end{equation}
\end{proposition}

\begin{proof}
For $s=-\frac{3}{4}+$ this was proved in Lemma 6.1 \cite{I-method}.
But it is easy to see that the arguments actually work for $s=-3/4$.
\end{proof}

Since $E_I^2(t)$ is very close to $E_I^4(t)$, then we will control
$E_I^4(t)$ and hence control $E_I^2(t)$. In order to control the
increase of $E_I^4(t)$, we need to control its derivative
\[\frac{d}{dt}E_I^4(t)=\Lambda_5(M_5),\]
where
\[M_5(\xi_1,\ldots,\xi_5)=-2i[\sigma_4(\xi_1,\xi_2,\xi_3,\xi_4+\xi_5)(\xi_4+\xi_5)]_{sym}.\]

\begin{proposition}
Assume $I\subset \R$ with $|I|\les 1$. Let $0\leq k_1\leq \ldots
\leq k_5$ and $k_4\geq 10$. Then we have
\begin{eqnarray}\label{eq:5linear1}
\left|\int_I\int \prod_{i=1}^5 P_{k_i}(w_i)(x,t)dxdt\right|\les
2^{\frac{5}{12}(k_1+k_2+k_3)}2^{-k_4}2^{-k_5}\prod_{j=1}^5\norm{\wh{P_{k_j}(w_j)}}_{X_{k_j}},
\end{eqnarray}
where if $k_j=0$ then $X_{k_j}$ is replaced by $\bar{X}_{k_j}$ on
the right-hand side.
\end{proposition}
\begin{proof}
From H\"older's inequality the left-hand side of \eqref{eq:5linear1}
is dominated by
\[\prod_{i=1}^3\norm{P_{k_i}(w_i)}_{L_x^3L_{t\in I}^\infty}\cdot \norm{P_{k_4}(w_4)}_{L_x^\infty L_{t}^2}\cdot \norm{P_{k_5}(w_5)}_{L_x^\infty L_{t}^2}.\]
For $\norm{P_{k_4}(w_4)}_{L_x^\infty L_{t}^2}$ and
$\norm{P_{k_5}(w_5)}_{L_x^\infty L_{t}^2}$ we use Proposition
\ref{propXkembedding}. For $\norm{P_{k_i}(w_i)}_{L_x^3L_{t\in
I}^\infty}$ we use interpolation between
$\norm{P_{k_i}(w_i)}_{L_x^2L_{t\in I}^\infty}$ and
$\norm{P_{k_i}(w_i)}_{L_x^4L_{t\in I}^\infty}$, and Proposition
\ref{propXkembedding}.
\end{proof}

\begin{proposition}
Let $\delta\les 1$. Assume $m$ is of the form \eqref{eq:m} with
$s=-3/4$, then
\begin{eqnarray}
\left|\int_0^\delta \Lambda_5(M_5;u_1,\ldots,u_5)dt\right|\les
N^{-\frac{15}{4}} \prod_{j=1}^5\norm{Iu_j}_{\bar F^{0}(\delta)}.
\end{eqnarray}
\end{proposition}
\begin{proof}
Following the proof of Lemma 5.2 \cite{I-method}, it suffices to
prove that
\begin{eqnarray*}
&&\sum_{k_1,\ldots, k_5\geq 0}\left|\int_0^\delta
\Lambda_5\left(\prod_{i=1}^3\frac{1}{(N+N_i)m(N_i)}\frac{1}{m(N_4)}\frac{1}{m(N_5)};P_{k_1}u_1,\ldots,P_{k_5}u_5\right)dt\right|\\
&&\les N^{-\frac{15}{4}}\prod_{i=1}^5
\norm{u_j}_{\bar{F}^0(\delta)}.
\end{eqnarray*}
Where $N_i=2^{k_i}$. From symmetry we may assume $N_1\geq N_2\geq
N_3$ and $N_4\geq N_5$ and two of the $N_i\ges N$. We fix the
extension $\wt{u}_i$ such that $\norm{\wt{u}_i}_{\bar{F}^0}\les
2\norm{{u}_i}_{\bar{F}^0(\delta)}$. For simplicity, we still denote
$u_i$.

The form \eqref{eq:m} with $s=-3/4$ implies that
$\frac{1}{(N+N_i)m(N_i)}\les N^{-3/4}\jb{N_i}^{-1/4}$ and
$\frac{1}{m(N_4)m(N_5)}\les N^{-3/2}N_4^{3/4}N_5^{3/4}$. Therefore
we need to control
\begin{eqnarray}\label{eq:deriE4}
N^{-\frac{15}{4}}\sum_{k_i}\int_0^\delta
\Lambda_5\left({\jb{N_1}^{-1/4}\jb{N_2}^{-1/4}\jb{N_3}^{-1/4}N_4^{3/4}N_5^{3/4}};u_1,\ldots,u_5\right)dt.
\end{eqnarray}
If $N_2\sim N_1 \ges N$, $N_4\les N_2$, we consider the worst case
$N_1\geq N_2\geq N_4\geq N_5\geq N_3$. From \eqref{eq:5linear1} we
get
\begin{eqnarray}
\eqref{eq:deriE4}&\les& N^{-\frac{15}{4}}\sum_{N_i}
\jb{N_1}^{-5/4}\jb{N_2}^{-5/4}\jb{N_3}^{1/6}N_4^{7/6}N_5^{7/6}\prod_{i=1}^5
\norm{\wh{P_{k_i}u}}_{X_{k_i}}\nonumber\\
&\les& N^{-\frac{15}{4}} \prod_{j=1}^5\norm{Iu_j}_{\bar
F^{0}(\delta)}.
\end{eqnarray}
The rest Cases $N_4\sim N_5 \ges N$, $N_1\les N_5$ or $N_1\sim N_4
\ges N$ follow in a similar ways. We omit the details.
\end{proof}

With these propositions, one can easily get global well-posedness of
the KdV equation using the same argument as in Section 6.4
\cite{I-method}. Moreover, we obtain that our global-in-time
solution satisfies
\begin{eqnarray}
\norm{u(t)}_{H^{-3/4}}\les \ (1+|t|)\cdot \norm{\phi}_{H^{-3/4}}.
\end{eqnarray}
 The proof for mKdV is just identical to the one in
\cite{I-method}, since it is easy to see that the Lemmas in Section
9.1 and 9.2 also hold for $s=1/4$.

\noindent{\bf Acknowledgment.} The author is very grateful to
Professor Carlos E. Kenig for encouraging the author to work on this
problem and helpful conversations, and to Professor Terence Tao for
the precious suggestions. The author also would like to thank
Professor Lizhong Peng and Professor Baoxiang Wang for the numerous
supports and encouragements. This work is supported in part by RFDP
of China No. 20060001010, the National Science Foundation of China,
grant 10571004; and the 973 Project Foundation of China, grant
2006CB805902, and the Innovation Group Foundation of NSFC, grant
10621061.

\end{document}